# FEEDBACK STABILIZATION OF A CLASS OF PERTURBED NONLINEAR AUTONOMOUS DIFFERENCE EQUATIONS

M. De la Sen
Institute of Research and Development of Processes. Campus of Leioa, Bizkaia, SPAIN
email: manuel.delasen@ehu.es

**Abstract**- This paper investigates the local asymptotic stabilization of a very general class of instable autonomous nonlinear difference equations which are subject to perturbed dynamics which can have a different order that that of the nominal difference equation. In the general case, the controller consists of two combined parts, namely, the feedback nominal controller which stabilizes the nominal (i.e. perturbation - free) difference equation plus an incremental controller which completes the stabilization in the presence of dynamics in the uncontrolled difference equation. A stabilization variant consists of using a single controller to stabilize the nominal difference equation and also the perturbed one under a smallness-type characterization of the perturbed dynamics. The study is based on Banach fixed point principle and it is also valid with slight modification for the stabilization of unstable oscillatory solutions.

## 1. Introduction

In this paper, the following nonautonomous difference equation is investigated:

$$\begin{aligned} x_n &= h_n\left(x_{n-1},\dots,x_{n-m}\right) \\ &= f_n\left(x_{n-1},\dots,x_{n-m_0}\right)+\tilde{f}_n\left(x_{n-1},\dots,x_{n-\tilde{m}}\right)+g_n\left(x_{n-1},\dots,x_{n-m_g}\right)+\tilde{g}_n\left(x_{n-1},\dots,x_{n-\tilde{m}_g}\right) \\ &= x_n^0+\tilde{x}_n+x_n^c+\tilde{x}_n^c\,;\quad n\in\mathbf{N} \end{aligned} \tag{1}$$

of order $m:=max\left(m_0,\tilde{m},m_g,\tilde{m}_g\right)\geq 1$ and initial conditions $x_{1-m},\dots,\ x_0$ where the four terms of the second identity are pair-wise identical in the same order as written, in which $h_n : D\subset\mathbf{R}^m\to\mathbf{R}$, $f_n : D_f\subset\mathbf{R}^{m_0}\to\mathbf{R}$, $\tilde{f}_n : D_{\tilde{f}}\subset\mathbf{R}^{\tilde{m}}\to\mathbf{R}$, $g_n : D_g\subset\mathbf{R}^{m_g}\to\mathbf{R}$ and $\tilde{g}_n : D_{\tilde{g}}\subset\mathbf{R}^{\tilde{m}_g}\to\mathbf{R}$ and $D$ is a nonempty subset of the union of the sets $D_f$, $D_{\tilde{f}}$, $D_g$ and $D_{\tilde{g}}$. The four pair-wise identical terms of the last identity have the following interpretations:

. $x_n^0$ is the nominal value of the uncontrolled nominal solution $x_n$ at the n-th sample in the absence of perturbations and controls

. $\tilde{x}_n$ is the perturbed uncontrolled solution which can be generated for perturbed parameterizations (then $\tilde{m}\leq m_0$) and possibly contributed by unmodeled dynamics (then $\tilde{m}>m_0$)

. $x_n^c$ is the correction by some nominal feedback controller of the uncontrolled nominal solution which can be potentially used to stabilize it or to improve it in some practical suitable sense provided it is already stable

. $\tilde{x}_n^c$ is the correction by adding some incremental feedback controller of the perturbed nominally controlled solution

. $\boldsymbol{N}$ is the set of natural numbers and $\boldsymbol{N}_0 = \boldsymbol{N} \cup \{0\}$ is that of nonnegative integers.

. $S^0$ and $cl\, S$ denote, respectively, the interior and closure of the set $S$ .

The stability and instability properties of nonautonomous difference equations has been investigated in [1-7], [17-18] and references there in. There is wide set of problems where stability of discrete systems involving either the discretization of time-continuous systems or being essentially digital by nature are of interest and involving very often the presence of nonlinearities. In those problems stability is commonly a required property for well-posedness. Among such problems, we can mention: a) those related to signal processing, b)models involving neural networks, c) adaptive control to deal with not perfectly known systems under combined estimation and control, d)problems related to modelling dynamic systems describing biological, medical or ecological systems, and e) those related to descriptions to rational difference equations. See, for instance, [9-15] and references there in. In particular, the nominal uncontrolled particular case of (1), namely, $f_n + \tilde{f}_n + g_n \equiv 0$ ; $\forall n \in \boldsymbol{N}$ has been studied in [1-4] under asymptotic stability conditions and dual instability versions of them, [1]. The objective of this manuscript is the generalization of the study of the stability and instability of equilibrium points of autonomous difference equations of [1-3] to the case of presence of additive perturbations (without formal distinction between parametrical perturbations, or structured or unstructured unmodeled dynamics). The perturbation- free difference equation will be referred to as the nominal uncontrolled one while the perturbed difference equation will be referred to as the uncontrolled perturbed difference equation. Two classes of feedback controllers are also proposed to stabilize the uncontrolled autonomous difference equation. The first class consists of two additive dynamics, namely, the nominal control for stabilization of the uncontrolled nominal equation plus an incremental controller for stabilization of the unmodeled dynamics. The second class consists of a single controller which stabilizes the whole uncontrolled dynamics for a certain tolerance to presence of perturbation dynamics of sufficiently small size characterized in terms of sufficiently small norm. The perturbed uncontrolled difference equation and the controlled difference equation can potentially possess distinct equilibrium points than the uncontrolled nominal difference equation. The formalism can also be applied to the study of feedback stabilization of unstable oscillations.

## 2. Vectorization preliminaries and linearization-based relations between equilibrium points and limit oscillatory solutions

Problems of major interest concerning (1) are (a) the characterization of a controller which stabilizes, at least locally around an equilibrium point, an unstable nominal difference equation, and (b) the stabilization of either a particular or a class of perturbed uncontrolled equations under a combined nominal plus incremental controller. It has to be pointed out that any equilibrium point of the uncontrolled equations can be re-allocated under a control action. In other words, the local stabilization via feedback control of an unstable equilibrium point of the uncontrolled equations may lead in parallel to a re-allocation of such an equilibrium point. An associate vector function to (1) is

$$V_h(u_1,....,u_m)=(h(u_1,....,u_m),u_1,...,u_{m-1})$$
$$=V_f(u_1,....,u_{m_0})+V_{\tilde{f}}(u_1,....,u_{\tilde{m}})+V_g(u_1,....,u_{m_g})+V_{\tilde{g}}(u_1,....,u_{\tilde{m}_g})$$
$$=\left(f(u_1,....,u_{m_0})+\tilde{f}(u_1,....,u_{\tilde{m}})+g(u_1,....,u_{m_g})+\tilde{g}(u_1,....,u_{\tilde{m}_g}),u_1,...,u_{m-1}\right) \quad (2)$$

, in particular $V_h(u_1)=h(u_1)$ if m=1, where

$$V_f(u_1,....,u_{m_0}):=\left(f(u_1,....,u_{m_0}),u_1,...,u_{m_0-1}\right);\ V_{\tilde{f}}(u_1,....,u_{\tilde{m}}):=\left(\tilde{f}(u_1,....,u_{\tilde{m}}),u_1,...,u_{\tilde{m}}\right)$$
$$V_g(u_1,....,u_{m_g}):=\left(g(u_1,....,u_{m_g}),u_1,...,u_{m_g-1}\right);\ V_{\tilde{g}}(u_1,....,u_{\tilde{m}_g})=\left(\tilde{g}(u_1,....,u_{m_{\tilde{g}}}),u_1,...,u_{m_{\tilde{g}}-1}\right) \quad (3)$$

The following result follows by simple inspection of (2) since $m\ge m_0$ :

**Lemma 2.1**. The vector function (2) can be expressed equivalently as

$$V_h(u_1,....,u_m)=$$
$$=\left(\bar{f}\left(u_1,....,u_{m_0},\overbrace{0,...,0}^{m-m_0}\right),u_1,...,u_{m-1}\right)+\left(\tilde{f}(u_1,....,u_{\tilde{m}})+g(u_1,....,u_{m_g})+\tilde{g}(u_1,....,u_{\tilde{m}_g}),\overbrace{0,...,0}^{m-1}\right) \quad (4)$$

where $\bar{f}:\left(D_f\cup D_{\tilde{f}}\right)\subset\boldsymbol{R}^m\to\boldsymbol{R}$ is defined from $f:D_f\subset\boldsymbol{R}^{m_0}\to\boldsymbol{R}$ by adding $(m-m_0)$ identically zero arguments. The set $D$ can be identical (although it is non-necessarily identical) to $D_f$ only if $m=m_0$ and then there is a unique such a mapping which is the identity self-mapping. □

The nominal and perturbed uncontrolled difference equations as well as the nominal controlled and perturbed controlled ones can have potentially distinct equilibrium points as follows. A generic " ad – hoc" description is also useful to describe some limit oscillatory solutions.

1) $\bar{x}^0$ is an equilibrium point of the uncontrolled nominal difference equation $x_n=f_n(x_{n-1},...,x_{n-m_0})$ if and only if $\bar{x}^0=f_n(\bar{x}^0,...,\bar{x}^0)$; $\forall n\in\boldsymbol{N}$. Then, $\bar{X}^0=(\bar{x}^0,...,\bar{x}^0)$ is the associate equilibrium point of the first-order autonomous $m_0$-order vector equation $X_n=V_{f_n}(X_{n-1})$; $\forall n\in\boldsymbol{N}$ obtained from the particular difference equation (1) $x_n=f_n(x_{n-1},...,x_{n-m_0})$ ; $\forall n\in\boldsymbol{N}$ via the nominal vector equation $V_{f_n}(u_1,....,u_m)=\left(f(u_1,....,u_{m_0}),u_1,....,u_{m_0}\right)$ provided that $V_{f_n}(D_f)\subseteq D_f$. A sequence solution $\left(\bar{x}_1^0,...\bar{x}_{m_0}^0\right)$ of $x_n=f_n(x_{n-1},...,x_{n-m_0})$ is a limit oscillatory solution of order at most $m_0$ if and only if $\bar{x}_{km_0+i}^0=f_n\left(\bar{x}_1^0,...\bar{x}_{m_0}^0\right)$; $\forall k,n\in\boldsymbol{N}$, $\forall i\in\overline{m}_0:=\{1,2,...,m_0\}$. Such a solution is trivially an equilibrium point if $\bar{x}_i^0=\bar{x}^0$; $\forall i\in\overline{m}_0$. The $m_0$ - real vector $\bar{X}^0=\left(\bar{x}_1^0,...\bar{x}_{m_0}^0\right)$ is the associate nominal limit oscillatory solution of order at most $m_0$ the first-order autonomous $m_0$-order vector equation $X_n = V_{f_n}(X_{n-1})$; $\forall n \in \boldsymbol{N}$ obtained from the particular difference equation (1) $x_n = f_n(x_{n-1},\dots,x_{n-m_0})$ ; $\forall n \in \boldsymbol{N}$ via the nominal vector equation $V_{f_n}(u_1,\dots,u_m) = \left(f(u_1,\dots,u_{m_0}), u_1,\dots,u_{m_0}\right)$.

2) $\bar{x}^{0p}$ is an equilibrium point of the uncontrolled perturbed difference equation $x_n = f_n(x_{n-1},\dots,x_{n-m_0}) + \tilde{f}_n(x_{n-1},\dots,x_{n-\tilde{m}})$ if and only if $\bar{x}^{0p} = f_n(\bar{x}^{0p},\dots,\bar{x}^{0p}) + \tilde{f}_n(\bar{x}^{0p},\dots,\bar{x}^{0p})$; $\forall n \in \boldsymbol{N}$ .Then, $\bar{X}^{0p} = \left(\bar{x}^{0p},\dots,\bar{x}^{0p}\right)$ is the associate equilibrium point of the first-order autonomous $m_{0p} := max(m_0, \tilde{m})$-order vector equation $X_n = V_{f_n + \tilde{f}_n}(X_{n-1})$; $\forall n \in \boldsymbol{N}$ according to

$$V_{f+\tilde{f}}\left(u_1,\dots,u_{m_{0p}}\right)$$

$$= \left( \bar{f}_0\left( u_1,\dots,u_{m_{0p}}, \overbrace{0,\ \dots,0}^{m_{0p}-m_0} \right), u_1,\dots,u_{m_{0p}-1} \right) + \left( \bar{\tilde{f}}_0\left( u_1,\dots,u_{\tilde{m}}, \overbrace{0,\ \dots,0}^{m_{0p}-\tilde{m}} \right), \overbrace{0,\ \dots,0}^{m_{0p}-1} \right) \qquad (5)$$

provided that $V_{f_n+\tilde{f}_n}\left(D_f \cup D_{\tilde{f}}\right) \subseteq D_f \cup D_{\tilde{f}}$ provided the set union is non-empty, where $\bar{f}_0$ , $\bar{\tilde{f}}_0 : \left(D_f \cup D_{\tilde{f}}\right) \subset \boldsymbol{R}^{m_{0p}} \to \boldsymbol{R}$ in (5), provided that $D_f \cup D_{\tilde{f}}$ is non-empty, take into account that the uncontrolled and nominal perturbed difference equations have potentially distinct orders and are built from $f : D_f \subset \boldsymbol{R}^{m_0} \to \boldsymbol{R}$, $\tilde{f} : D_{\tilde{f}} \subset \boldsymbol{R}^{\tilde{m}} \to \boldsymbol{R}$ as in the parallel construction of Lemma 2.1, Eq.(4). We can describe limit oscillatory solutions of the uncontrolled perturbed difference equation of order at most $m_{0p}$, or equivalently those of its associate vector function, by a sequence solution $\bar{x}^{0p}_{km_{0p}+i} = f_n\left(\bar{x}^{0p}_1,\dots\bar{x}^{0p}_{m_{0p}}\right) + \tilde{f}_n\left(\bar{x}^{0p}_1,\dots\bar{x}^{0p}_{m_{0p}}\right)$; $\forall k,n \in \boldsymbol{N}$, $\forall i \in \overline{m}_{0p}$.

3) $\bar{x}^c$ is an equilibrium point of the controlled nominal difference equation $x_n = f_n(x_{n-1},\dots,x_{n-m_0}) + g_n(x_{n-1},\dots,x_{n-m_g})$ if and only if $\bar{x}^c = f_n(\bar{x}^c,\dots,\bar{x}^c) + g_n(\bar{x}^c,\dots,\bar{x}^c)$; $\forall n \in \boldsymbol{N}$ . Then, $\bar{X}^c = \left(\bar{x}^c,\dots,\bar{x}^c\right)$ is the associate equilibrium point of the first-order autonomous $m_c := max(m_0, m_g)$-order vector equation $X_n = V_{f_n+g_n}(X_{n-1})$; $\forall n \in \boldsymbol{N}$ , provided that $V_{f+g}\left(D_f \cup D_g\right) \subseteq D_f \cup D_g$ , provided that such a union is non-empty. The equivalent first-order vector equations are defined via the associate vector function $V_{f+g}\left(D_f \cup D_g\right)$ through ad-hoc functions $\bar{f}_c : D \to D_f$ and $\bar{g}_c : D \to D_g$ built according to the corresponding associate vector equation defined in a similar way to (4) and (5). This is directly extended to limit oscillatory solutions of order at most $m_c$ of the controlled nominal difference equation, which can be equivalently expressed in vector form, in the same way as above.

4) $\bar{x}^{cp}$ is an equilibrium point of the controlled perturbed difference equation (1) if and only if

$$\bar{x}^{cp} = h_n\left(\bar{x}^{cp},\dots,\bar{x}^{cp}\right)$$
$$= f_n\left(\bar{x}^{cp},\dots,\bar{x}^{cp}\right) + \tilde{f}_n\left(\bar{x}^{cp},\dots,\bar{x}^{cp}\right) + g_n\left(\bar{x}^{cp},\dots,\bar{x}^{cp}\right) + \tilde{g}_n\left(\bar{x}^{cp},\dots,\bar{x}^{cp}\right);\ \forall n \in \boldsymbol{N} \qquad (6)$$

Then, $\bar{X}^{cp} = \left(\bar{x}^{cp},\dots,\bar{x}^{cp}\right)$ is the associate equilibrium point of the first-order autonomous m-order vector equation $X_n = V_{h_n}\left(X_{n-1}\right)$; $\forall n \in \boldsymbol{N}$ defined via (4) provided that $V_{h_n}(D) \subseteq D$. A sequence solution $\left(\bar{x}_1^{cp},\dots\bar{x}_m^{cp}\right)$ of $x_n = f_n\left(x_{n-1},\dots,x_{n-m}\right)$ is a limit oscillatory solution of order at most $m$ if and only if $\bar{x}_{km+i}^{cp} = f_n\left(\bar{x}_1^{cp},\dots\bar{x}_m^{cp}\right)$; $\forall k,n \in \boldsymbol{N}$; $\forall i \in \bar{m}$. Such a solution is trivially an equilibrium point if $\bar{x}_i^{cp} = \bar{x}^{cp}$; $\forall i \in \bar{m}$. The $m$ - real vector $\bar{X}^{cp} = \left(\bar{x}_1^{cp},\dots\bar{x}_m^{cp}\right)$ is the associate nominal limit oscillatory solution of order at most $m$ the first-order autonomous $m$-order vector equation $X_n = V_{f_n}\left(X_{n-1}\right)$; $\forall n \in \boldsymbol{N}$ obtained from the particular difference equation (1) $x_n = f_n\left(x_{n-1},\dots,x_{n-m}\right)$; $\forall n \in \boldsymbol{N}$ via the nominal vector equation $V_{f_n}\left(u_1,\dots,u_m\right) = \left(f\left(u_1,\dots,u_m\right),u_1,\dots,u_m\right)$.

*Remark 2.1.* Note that the above description allows the characterization of equilibrium points as particular cases of limit oscillatory solutions. Note also that limiting oscillatory solutions can exceed the order of the difference equations if such a solution has a repeated pattern of more elements than the order of the difference equations. Details are omitted since the analysis methods is close to the above one in both scalar and equivalent vector forms. Limiting oscillatory solutions are relevant in some applications, in particular, in the fields of Communications, design of electronic oscillators, etc. □

*Remark 2.2*: The difference equation $x_n = f_n\left(x_{n-1},\dots,x_{n-m}\right)$; $\forall n \in \boldsymbol{N}$ has been pointed to be equivalent to its associate vector equation $X_n = V_{f_n}\left(X_{n-1}\right)$; $\forall n \in \boldsymbol{N}$. Then, the nominal uncontrolled difference equation admits the representation $x_n = f_n\left(X_{n-1}\right) = f_n\left(V_{f_{n-1}}\left(X_{n-2}\right)\right)$; $\forall n \in \boldsymbol{N}$. Proceeding recursively:

$$x_n = f_n\left(X_{n-1}\right) = f_n\left(V_{f_{n-1}}\left(X_{n-2}\right)\right) = f_n\left(G_{n-1}^{f}\left(X_0\right)\right);\ \forall n \in \boldsymbol{N}$$

By defining $G_n^f := V_{f_n} \mathrm{o} V_{f_{n-1}} \mathrm{o} \dots \mathrm{o} V_{f_1}$; $\forall n \in \boldsymbol{N}_0$ with $G_0^f$ being identity, [1]. Close composed mappings to describe the various uncontrolled and controlled (nominal or ) versions of (1) are:

$$G_0^{f+\tilde{f}} = \mathrm{id},\ G_n^{f+\tilde{f}} := V_{f_n+\tilde{f}_n} \mathrm{o} V_{f_{n-1}+\tilde{f}_{n-1}} \mathrm{o} \dots \mathrm{o} V_{f_1+\tilde{f}_1};\ \forall n \in \boldsymbol{N}$$

$$G_0^{f+\tilde{f}+g} = \mathrm{id},\ G_n^{f+\tilde{f}+g} := V_{f_n+\tilde{f}_n+g_n} \mathrm{o} V_{f_{n-1}+\tilde{f}_{n-1}+g_{n-1}} \mathrm{o} \dots \mathrm{o} V_{f_1+\tilde{f}_1+g_1};\ \forall n \in \boldsymbol{N}$$

$$G_0^{h} = \mathrm{id},\ G_n^{h} := V_{h_n} \mathrm{o} V_{h_{n-1}} \mathrm{o} \dots \mathrm{o} V_{h_1};\ \forall n \in \boldsymbol{N}$$ □

The following result is direct by inspection of (1):

**Proposition 2.2**. The nominal and perturbed uncontrolled associate vector functions may have a common equilibrium point or a common limiting oscillation of order at most m only if $m = m_0 = max\left(m_0, \tilde{m}\right)$. The nominal and perturbed uncontrolled associate vector functions as well as the controlled and perturbed controlled ones may have a common equilibrium point only if, in addition, $m_0 = max\left(m_0, m_g, \tilde{m}_g\right) = max\left(m_0, \tilde{m}_g\right)$.

**Proof**: If the conditions fail and the vector functions referred to have some common equilibrium point, this one, should have different dimension depending on the equation what is a contradiction. The proof is also valid "mutatis-mutandis" for limiting oscillation of orders at most *m*. □

It is now discussed the presence of limit oscillations of the uncontrolled perturbed difference equations in a neighbourhood centred about a nominal limit oscillation. A similar analysis is also useful for closeness of the limiting oscillatory solutions to that of a given difference equation of any of the three remaining difference ones under investigation.

**Theorem 2.3**. Assume the following:

1) $m_p := max\left(m_0, \tilde{m}\right) = m_0$, $V_{f_n}\left(D_f\right) \subseteq D_f \neq \varnothing$ and $V_{f_n + \tilde{f}_n}\left(D_f \cup D_{\tilde{f}}\right) \subseteq D_f \cup D_{\tilde{f}} \neq \varnothing$; $\forall n \in \boldsymbol{N}$ where $V_{f_n + \tilde{f}_n}\left(u_1, ..., u_{m_p}\right)$ is defined in (5); $\forall n \in \boldsymbol{N}$.

2) $\left.\frac{\partial f_n(X)}{\partial X^T}\right|_{\overline{X}^0}$, $\left.\frac{\partial \tilde{f}_n(X)}{\partial X^T}\right|_{\overline{X}^0}$ exist within a neighborhod of $\overline{X}^0$, which is a limit oscillatory solution of the vector uncontrolled nominal equation of order at most $m_p = m_0$ (including potentially nominal equilibrium points as particular cases).

3) The inverse $m_0$ -matrix $\left(I_{m_p} - M^0\left(\overline{X}^0\right)\right)^{-1}$ exists, where $M^0\left(\overline{X}^0\right) = \left.\frac{\partial V_f(X)}{\partial X^T}\right|_{\overline{X}^0} + \left.\frac{\partial V_{\tilde{f}}(X)}{\partial X^T}\right|_{\overline{X}^0}$, and $I_{m_0}$ is the $m_0$-identity matrix of $\boldsymbol{R}^{m_0}$.

Then $\hat{\overline{X}}^{0p} = \overline{X}^0 + \left(I_m - M^0\left(\overline{X}^0\right)\right)^{-1} V_{\tilde{f}}\left(\overline{X}^0\right)$ is a linear estimate of limit oscillatory solutions of order at most $m_0$ (including, as particular cases, potential equilibrium points) of the vector perturbed difference equation. The estimate is closes to true values of $\overline{X}^{0p}$ as $\left\|V_{\tilde{f}}\left(\overline{X}^0\right)\right\|$ is sufficiently small.

If $rank\left(I_{m_p} - M^0\left(\overline{X}^0\right), V_{\tilde{f}}\left(\overline{X}^0\right)\right) = rank\left(I_{m_p} - M^0\left(\overline{X}^0\right)\right) < m_p$ then there are infinitely many first-order estimates $\hat{\overline{X}}^{0p}$ of limiting oscillatory solutions of the vector uncontrolled nominal equation of order at most $m_0$. If $m_p - 1 \geq rank\left(I_m - M^0\left(\overline{X}^0\right), V_{\tilde{f}}\left(\overline{X}^0\right)\right) > rank\left(I_{m_p} - M^0\left(\overline{X}^0\right)\right)$ then there is no such an estimate.

**Proof**: Note that $m_p = m_0$ implies that $\tilde{m} \le m_0$. Define $\Delta\overline{X}^{0p,0} := \overline{X}^{0p} - \overline{X}^{0}$ which is rewritten below after using a linearized perturbed difference vector equation since the perturbed equilibrium point is within a neighbourhood of the nominal one

$$\Delta\overline{X}_{n+1}^{0p,0} = M_n^0\left(\overline{X}^0\right)\Delta\overline{X}_n^{0p,0} + V_{\tilde{f}_n}\left(\overline{X}^0\right) + o\left(\left\|\Delta\overline{X}_n^{0p,0}\right\|\right)I_{m_p}$$

$$= M^0\left(\overline{X}^0\right)\Delta\overline{X}_n^{0p,0} + V_{\tilde{f}}\left(\overline{X}^0\right) + o\left(\left\|\Delta\overline{X}_n^{0p,0}\right\|\right)I_{m_p} \; ; \; \forall n \in \boldsymbol{N} \tag{7}$$

where

$$M_n^0\left(\overline{X}^0\right) = \left.\frac{\partial V_{f_n}(X)}{\partial X^T}\right|_{\overline{X}^0} + \left.\frac{\partial V_{\tilde{f}_n}(X)}{\partial X^T}\right|_{\overline{X}^0} = \begin{bmatrix} \left.\frac{\partial f_n(X)}{\partial X^T}\right|_{\overline{X}^0} + \left.\frac{\partial \tilde{f}_n(X)}{\partial X^T}\right|_{\overline{X}^0} \\ I_{m_p-1} \quad 0 \end{bmatrix}$$

$$= M^0\left(\overline{X}^0\right) = \left.\frac{\partial V_f(X)}{\partial X^T}\right|_{\overline{X}^0} + \left.\frac{\partial V_{\tilde{f}}(X)}{\partial X^T}\right|_{\overline{X}^0} = \begin{bmatrix} \left.\frac{\partial f(X)}{\partial X^T}\right|_{\overline{X}^0} + \left.\frac{\partial \tilde{f}(X)}{\partial X^T}\right|_{\overline{X}^0} \\ I_{m-1} \quad 0 \end{bmatrix} ; \; \forall n \in \boldsymbol{N} \tag{8}$$

if $m_0 \ge 2$, where $I_{m_0-1}$ is the $(m_0-1)$ identity matrix and superscript $T$ denotes transposition, since $V_{f_n}\left(\overline{X}^0\right) = V_f\left(\overline{X}^0\right)$, $V_{\tilde{f}_n}\left(\overline{X}^0\right) = V_{\tilde{f}}\left(\overline{X}^0\right)$, $f_n\left(\overline{X}^0\right) = f\left(\overline{X}^0\right)$, $\tilde{f}_n\left(\overline{X}^0\right) = \tilde{f}\left(\overline{X}^0\right)$, $\forall n \in \boldsymbol{N}$ ; and

$$M_n^0\left(\overline{X}^0\right) = \left.\frac{\partial V_{f_n}(X)}{\partial X}\right|_{\overline{X}^0} + \left.\frac{\partial V_{\tilde{f}_n}(X)}{\partial X}\right|_{\overline{X}^0} = M^0\left(\overline{X}^0\right) = \left.\frac{\partial V_f(X)}{\partial X}\right|_{\overline{X}^0} + \left.\frac{\partial V_{\tilde{f}}(X)}{\partial X}\right|_{\overline{X}^0} ; \; \forall n \in \boldsymbol{N} \tag{9}$$

if $m_0 = 1$. Taking $\Delta\overline{X}_{n+1}^{0p,0} = \Delta\overline{X}_n^{0p,0} = \Delta\overline{X}^{0p,0}$ ; $\forall n \in \boldsymbol{N}$, one gets from (7):

$$\Delta\overline{X}^{0p,0} = \left(I_{m_0} - M^0\left(\overline{X}^0\right)\right)^{-1} V_{\tilde{f}}\left(\overline{X}_0\right) + o\left(\left\|\Delta\overline{X}^{0p,0}\right\|\right)I_{m_0} \tag{10}$$

provided that $\left(I_{m_0} - M^0\left(\overline{X}^0\right)\right)^{-1}$ exists so that $\Delta\hat{\overline{X}}^{0p,0} = \left(I_m - M^0\left(\overline{X}^0\right)\right)^{-1} V_{\tilde{f}}\left(\overline{X}_0\right)$ is an estimate of $\Delta\overline{X}_{n+1}^{0p,0}$ so that if $\left\|M^0\left(\overline{X}^0\right)\right\| < 1$ and $\left\|V_{\tilde{f}}\left(\overline{X}_0\right)\right\| < \left(1 - \left\|M^0\left(\overline{X}^0\right)\right\|\right)\varepsilon$ for some $\varepsilon \in \boldsymbol{R}_+$, then one gets from Banach´s perturbation lemma,[16]:

$$\left\|\Delta\hat{\overline{X}}^{0p,0}\right\| \le \left\|\left(I_m - M^0\left(\overline{X}^0\right)\right)^{-1}\right\| \left\|V_{\tilde{f}}\left(\overline{X}_0\right)\right\| \le \frac{\left\|V_{\tilde{f}}\left(\overline{X}^0\right)\right\|}{1 - \left\|M^0\left(\overline{X}^0\right)\right\|} < \varepsilon \tag{11}$$

Since for $\varepsilon = 0$, $V_{\tilde{f}}\left(\overline{X}_0\right) = \Delta\overline{X}^{0p,0} = 0$ then, for a sufficiently small $\varepsilon^*$ such that $\frac{\left\|V_{\tilde{f}}\left(\overline{X}^0\right)\right\|}{1 - \left\|M^0\left(\overline{X}^0\right)\right\|} < \varepsilon^*$ and for any $\varepsilon \le \varepsilon^*$, $o\left(\left\|\Delta\overline{X}^{0p,0}\right\|\right) \le \varepsilon/2$ what occurs in particular, for $\varepsilon < 1$ if $f, \tilde{f} : \boldsymbol{R}^m \to \boldsymbol{R}$ are furthermore analytic in an open ball of $\boldsymbol{R}^m$ centred at $\overline{X}_0$ of radius $\rho = 3\varepsilon/2$. The conditions for the existence of infinitely many first-order estimates of $\Delta\hat{\overline{X}}^{0p,0}$ or the existence of none of them is direct from compatible and incompatible conditions for linear algebraic systems of equations according to Rouché – Froebenius theorem from Linear Algebra. □

Note that it can occur for the nominal and perturbed uncontrolled difference equations to have common equilibrium points. On the other hand, it is possible to obtain linear similar first-order comparison results to those of Theorem 2.3 for the estimates of the equilibrium points of the corrected closed-loop system via an incremental controller related to those of the controlled system without incremental controller. An "ad-hoc" result is now stated without proof which can be performed very closely to that of Theorem 2.3:

**Theorem 2.4**. Assume the following:

1) $m = max(m_0, \tilde{m}, m_g)$ and $\tilde{m}_g \le max(m_0, \tilde{m}, m_g)$, so that $m = max(m_0, \tilde{m}, m_g, \tilde{m}_g)$, and $V_{h_n}(D) \subseteq D \neq \varnothing$, $V_{f_n+\tilde{f}_n+g_n}\left(D_f \cup D_{\tilde{f}} \cup D_g\right) \subseteq D_f \cup D_{\tilde{f}} \cup D_g \neq \varnothing$; $\forall n \in \boldsymbol{N}$; where $V_{f_n+\tilde{f}_n+g_n}(u_1, \ldots, u_m)$ is defined correspondingly to (5) for this case; $\forall n \in \boldsymbol{N}$.

2) $\left.\dfrac{\partial\left(f_n(X) + \tilde{f}_n(X) + g_n(X)\right)}{\partial X^T}\right|_{\overline{X}^c}$, $\left.\dfrac{\partial \tilde{g}_n(X)}{\partial X^T}\right|_{\overline{X}^c}$ exist within a neighborhod of $\overline{X}^c$, which is a limit oscillatory solution of order at most $m$ of the vector controlled nominal equation, i.e. the vector nominal uncontrolled via feedback of the nominal controller (including potentially nominal equilibrium points).

3) The inverse m-matrix $\left(I_m - M^c\left(\overline{X}^c\right)\right)^{-1}$ exists, where $M^c\left(\overline{X}^c\right) = \left.\dfrac{\partial V_h(X)}{\partial X^T}\right|_{\overline{X}^c}$.

Then $\hat{\overline{X}}^{cp} = \overline{X}^c + \left(I_m - M^c\left(\overline{X}^c\right)\right)^{-1} V_{\tilde{g}}\left(\overline{X}^c\right)$ is a linear estimate of limit oscillatory solutions of order at most $m$ (including potential equilibrium points as particular cases) of the vector controlled difference equation under the combined nominal and correction controllers from its corresponding counterpart under the nominal controller only. The estimate closes to true values as $\left\| V_{\tilde{g}}\left(\overline{X}^c\right)\right\|$ is sufficiently small.

If $rank\left(I_m - M^c\left(\overline{X}^c\right), V_{\tilde{g}}\left(\overline{X}^c\right)\right) = rank\left(I_m - M^c\left(\overline{X}^c\right)\right) < m$ then there are infinitely many first-order estimates $\hat{\overline{X}}^{cp}$ of limiting oscillatory solutions of the vector uncontrolled nominal equation of order at most $m$. If $m - 1 \ge rank\left(I_m - M^c\left(\overline{X}^c\right), V_{\tilde{f}}\left(\overline{X}^c\right)\right) > rank\left(I_m - M^c\left(\overline{X}^c\right)\right)$ then there is no such an estimate. □

The same linearization technique can be used to compare closely allocated equilibrium points of the same dimension for other pairs of the involved systems. In this way, the following results follow, respectively, for the nominal uncontrolled and controlled difference equations and for the uncontrolled nominal and controlled perturbed ones and, equivalently, for the associate pairs of vector systems as follows:

**Theorem 2.5**. Assume the following:

1) $m_c := max(m_0, m_g) = m_0$, $V_{f_n}\left(D_f\right) \subseteq D_f \neq \varnothing$ and $V_{f_n+g_n}\left(D_f \cup D_g\right) \subseteq D_f \cup D_g \neq \varnothing$; $\forall n \in \boldsymbol{N}$ where $V_{f_n+g_n}\left(u_1, \ldots, u_{m_c}\right)$ is defined correspondingly to (5) for this case; $\forall n \in \boldsymbol{N}$.

2) $\left.\frac{\partial f_n(X)}{\partial X^T}\right|_{\overline{X}^0}$, $\left.\frac{\partial g_n(X)}{\partial X^T}\right|_{\overline{X}^0}$ exist within a neighborhod of $\overline{X}^0$, which is a limit oscillatory solution of the vector uncontrolled nominal equation of order at most $m_0$ (including potentially nominal equilibrium points as particular cases).

3) The inverse $m_0$ -matrix $\left(I_{m_0} - M^c\left(\overline{X}^0\right)\right)^{-1}$ exists, where $M^c\left(\overline{X}^0\right) = \left.\frac{\partial V_f(X)}{\partial X^T}\right|_{\overline{X}^0} + \left.\frac{\partial V_g(X)}{\partial X^T}\right|_{\overline{X}^0}$,

and $I_{m_0}$ is the $m_0$-identity matrix of $\boldsymbol{R}^{m_0}$.

Then $\hat{\overline{X}}^c = \overline{X}^0 + \left(I_m - M^c\left(\overline{X}^0\right)\right)^{-1} V_g\left(\overline{X}^0\right)$ is a linear estimate of limit oscillatory solutions of order at most $m_0$ (including, as particular cases, potential equilibrium points) of the vector controlled difference equation from its nominal uncontrolled counterpart. The estimate closes to true values as $\left\| V_g\left(\overline{X}^0\right)\right\|$ is sufficiently small.

If $rank\left(I_{m_0} - M^c\left(\overline{X}^0\right), V_g\left(\overline{X}^0\right)\right) = rank\left(I_{m_0} - M^c\left(\overline{X}^0\right)\right) < m_0$ then there are infinitely many first-order estimates $\hat{\overline{X}}^c$ of limiting oscillatory solutions of the vector controlled nominal equation of order at most $m_c = m_0$. If $m_0 - 1 \geq rank\left(I_{m_0} - M^c\left(\overline{X}^0\right), V_g\left(\overline{X}^0\right)\right) > rank\left(I_{m_0} - M^c\left(\overline{X}^0\right)\right)$ then there is no such an estimate. □

**Theorem 2.6**. Assume the following:

1) $m_{cp} := max\left(m_0, \tilde{m}_0, m_g\right) = max\left(m_0, \tilde{m}_0\right) = m_0$, $V_{f_n}\left(D_f \cup D_{\tilde{f}}\right) \subseteq D_f \cup D_{\tilde{f}} \neq \varnothing$, and $V_{f_n + \tilde{f}_n + g_n}\left(D_f \cup D_{\tilde{f}} \cup D_g\right) \subseteq D_f \cup D_{\tilde{f}} \cup D_g \neq \varnothing$; $\forall n \in \boldsymbol{N}$
where $V_{f_n + \tilde{f}_n + g_n}\left(u_1, ..., u_{m_c}\right)$ is defined correspondingly to (5) for this case; $\forall n \in \boldsymbol{N}$.

2) $\left.\frac{\partial f_n(X)}{\partial X^T}\right|_{X^0}$, $\left.\frac{\partial \tilde{f}_n(X)}{\partial X^T}\right|_{X^0}$, $\left.\frac{\partial g_n(X)}{\partial X^T}\right|_{X^0}$ exist within a neighborhod of $X^0$, which is a limit oscillatory solution of the vector uncontrolled nominal equation of order at most $m_0$ (including potentially nominal equilibrium points as particular cases).

3) The inverse $m_0$ -matrix $\left(I_{m_0} - M^{cp}\left(X^0\right)\right)^{-1}$ exists, where $M^{cp}\left(X^0\right) = \left.\frac{\partial V_f(X)}{\partial X^T}\right|_{X^0} + \left.\frac{\partial V_{\tilde{f}}(X)}{\partial X^T}\right|_{X^0} + \left.\frac{\partial V_g(X)}{\partial X^T}\right|_{X^0}$, and $I_{m_0}$ is the $m_0$-identity matrix of $\boldsymbol{R}^{m_0}$.

Then $\hat{\overline{X}}^{cp} = X^0 + \left(I_m - M^{cp}\left(X^0\right)\right)^{-1} V_{\tilde{f}+g}\left(X^0\right)$ is a linear estimate of limit oscillatory solutions of order at most $m_0$ (including, as particular cases, potential equilibrium points) of the vector controlled difference equation from its nominal uncontrolled counterpart. The estimate closes to true values as $\left\| V_{\tilde{f}+g}\left(X^0\right)\right\|$ is sufficiently small.

If $rank\left(I_{m_0}-M^{cp}\left(X^{0}\right),V_{\tilde{f}+g}\left(X^{0}\right)\right)=rank\left(I_{m_0}-M^{cp}\left(X^{0}\right)\right)<m_0$ then there are infinitely many first-order estimates $\hat{\overline{X}}^{cp}$ of limiting oscillatory solutions of the vector controlled nominal equation of order at most $m_{cp}=m_0$. If $m_0-1\ge rank\left(I_{m_0}-M^{cp}\left(X^{0}\right),V_{\tilde{f}+g}\left(X^{0}\right)\right)>rank\left(I_{m_0}-M^{cp}\left(X^{0}\right)\right)$ then there is no such an estimate. □

**3. Some stability and instability properties**

The following result holds concerning the stabilization via a feedback controller of an unstable uncontrolled equilibrium point. The controller consists, in general , of two parts, namely a) the nominal controller used to stabilize the uncontrolled difference equation; and b) the incremental controller used to stabilize the difference equation which includes perturbed parameters and/or perturbed dynamics. The stabilization process admits the double interpretation of the above section in terms of stabilization of either equilibrium points or that of oscillatory solutions. The equilibrium points can potentially vary under perturbations and the presence of feedback controllers.

**Theorem 3.1**. The following properties hold:

**(i)** Let $\overline{x}^{0}$ and $\overline{x}^{c}$ be two equilibrium points of the nominal and nominal controlled difference equations with corresponding ones $\overline{X}^{0}$ and $\overline{X}^{c}$ in the associate vector equations. Assume that $m=m_0=m_c$, $\overline{X}^{0}\in cl\,S\cap A_0$ and $\overline{X}^{c}\in cl\,S\cap A_p$ with $\varnothing\neq S^0\subseteq S\subseteq A_0\cap A_c$, where $S$ is an invariant subset of solutions of the associate vector equations for all $V_{f_n}$ and $V_{f_n+g_n}$; i.e. $V_{f_n}(S)\subseteq S$ and $V_{f_n+g_n}(S)\subseteq S$, for any controller in $\boldsymbol{C}$ consisting in a nominal controller, with

$$A_0:=\left\{X\in R^{m}:\left|f_n(X)-\overline{x}^{0}\right|\ge\alpha_n\left\|X-\overline{X}^{0}\right\|;\forall n\in\boldsymbol{N}\right\}\tag{12}$$

$$A_c:=\left\{X\in R^{m}:\left|f_n(X)+g_n(X)-\overline{x}^{c}\right|\le\beta_n^{c}\left\|X-\overline{X}^{c}\right\|;\forall n\in\boldsymbol{N}\right\}\tag{13}$$

for some real nonnegative sequences $\{\alpha_i\}_{i\in\boldsymbol{N}}$ and $\{\beta_i^{c}\}_{i\in\boldsymbol{N}}$. If $\{\overline{\alpha}_n\}_{n\in\boldsymbol{N}}$ is unbounded, where $\overline{\alpha}_n=\prod_{i=1}^{n}\alpha_i$, and $\{\overline{\beta}_n^{c}\}_{n\in\boldsymbol{N}}$ is such that $\limsup_{n\to\infty}\overline{\beta}_n^{c}<1$, where $\overline{\beta}_n^{c}=\prod_{i=1}^{n}\beta_i^{c}$, then $\overline{X}^{0}$ is unstable where $\overline{X}^{c}$ is locally asymptotically stable with respect to S.

**(ii)** Let $\overline{x}^{0p}$ and $\overline{x}^{c}$ be two equilibrium points of the uncontrolled perturbed and nominal controlled (via the nominal plus the incremental controllers) difference equations with corresponding ones $\overline{X}^{0p}$ and $\overline{X}^{c}$ in the associate vector equations. If $m=m_{0p}=max(m_0,\tilde{m})=m_{cp}$, $\overline{X}^{0p}\in cl\,S\cap A_{0p}$ and $\overline{X}^{c}\in cl\,S\cap A_c$ with $\varnothing\neq S^{0}\subseteq S\subseteq A_{0p}\cap A_c$, where $S$ is an invariant subset of the solutions of the associate vector equations for all $V_{f_n+\tilde{f}_n}$ and $V_{f_n+\tilde{f}_n+g_n}$; i.e. $V_{f_n+\tilde{f}_n}(S)\subseteq S$ and $V_{f_n+\tilde{f}_n+g_n}(S)\subseteq S$, with

$$A_{0p} := \left\{ X \in R^{m} : \left| f_n(X) + \tilde{f}_n(X) - x^{0\bar{p}} \right| \geq \alpha_{pn} \left\| X - \bar{X}^{0p} \right\| ; \forall n \in \boldsymbol{N} \right\} \tag{14}$$

$$A_{c} := \left\{ X \in R^{m} : \left| f_n(X) + \tilde{f}_n(X) + g_n(X) - \bar{x}^{c} \right| \leq \beta_n^c \left\| X - \bar{X}^{c} \right\| ; \forall n \in \boldsymbol{N} \right\} \tag{15}$$

for some nonnegative sequences $\{\alpha_{pi}\}_{i\in \boldsymbol{N}}$ and $\{\beta_i^c\}_{i\in \boldsymbol{N}}$. Thus, if $\{\bar{\alpha}_{pn}\}_{n\in \boldsymbol{N}}$ is unbounded, where $\bar{\alpha}_{pn} = \prod_{i=1}^{n} \alpha_{pi}$, and $\{\bar{\beta}_n^c\}_{n\in \boldsymbol{N}}$ is such that $\limsup_{n\to\infty} \bar{\beta}_n^c < 1$, where $\bar{\beta}_n^c = \prod_{i=1}^{n} \beta_i^c$, then $\bar{X}^{0p}$ is unstable where $\bar{X}^{c}$ is locally asymptotically stable with respect to *S*.

**(iii)** Let $\bar{x}^{0p}$ and $\bar{x}^{cp}$ be two equilibrium points of the perturbed and perturbed controlled (via the nominal plus the incremental controllers) difference equations with corresponding ones $\bar{X}^{0p}$ and $\bar{X}^{cp}$ in the associate vector equations. If $m = m_{0p} = m_{cp}$, $\bar{X}^{0p} \in cl\, S \cap A_{0p}$ and $\bar{X}^{cp} \in cl\, S \cap A_{cp}$ with $\emptyset \neq S^{0} \subseteq S \subseteq A_{0p} \cap A_{cp}$, where $S$ is an invariant subset of the solutions of the associate vector equations for all $V_{f_n+\tilde{f}_n}$ and $V_{h_n}$, i.e. $V_{f_n+\tilde{f}_n}(S) \subseteq S$ and $V_{h_n}(S) \subseteq S$, for some nonnegative real sequences $\{\alpha_i\}_{i\in \boldsymbol{N}}$ and $\{\beta_{pi}^c\}_{i\in \boldsymbol{N}}$, and

$$A_{cp} := \left\{ X \in R^{m} : \left| h_n(X) - \bar{x}^{cp} \right| \leq \beta_{pn}^c \left\| X - \bar{X}^{cp} \right\| ; \forall n \in \boldsymbol{N} \right\} \tag{16}$$

Thus, if $\{\bar{\alpha}_n\}_{n\in \boldsymbol{N}}$ is unbounded, where $\bar{\alpha}_n = \prod_{i=1}^{n} \alpha_i$, and $\{\bar{\beta}_{pn}^c\}_{n\in \boldsymbol{N}}$ is such that $\limsup_{n\to\infty} \bar{\beta}_{pn}^c < 1$, where $\bar{\beta}_{pn}^c = \prod_{i=1}^{n} \beta_{pi}^c$, then $\bar{X}^{0p}$ is unstable where $\bar{X}^{cp}$ is locally asymptotically stable with respect to *S*. □

The above result is directly extendable to stabilization of unstable oscillatory solutions to the light of the former discussions in Section 2. Explicit conditions for the fulfilment of Theorem 3.1(iii), which imply the local asymptotic stabilization within an invariant set around the equilibrium points of the unstable perturbed uncontrolled system, are given in the subsequent result. The stabilization mechanism is achieved by synthesizing a controller consisting of combined nominal controller with an incremental controller. The nominal controller stabilized the nominal difference equation in the absence of perturbations while the incremental one completes the stabilization for the whole uncontrolled difference equation.

**Theorem 3.2**. Assume that $m = m_{0p} = m_{cp}$ with $\bar{X}^{0p} \in cl\, S \cap A_{0}$ and $\bar{X}^{cp} \in cl\, S \cap A_{cp}$ being unique equilibrium points in $cl\, S \cap A_{0}$, respectively in $cl\, S \cap A_{cp}$, where $S \subseteq A_{0p} \cap A_{cp}$ is invariant under all $V_{f_n+\tilde{f}_n}$ and $V_{h_n}$ for a class of controllers $\boldsymbol{C}$; i.e. $V_{f_n+\tilde{f}_n}(S) \subseteq S$ and $V_{h_n}(S) \subseteq S$ for any combined nominal plus incremental controller in the class $\boldsymbol{C}$; $\forall n \in \boldsymbol{N}$. Define $\Delta \bar{x}^{0p} := \bar{x}^{cp} - \bar{x}^{0p}$ and $\Delta \bar{X}^{0p} := \bar{X}^{cp} - \bar{X}^{0p}$ being sufficiently close to zero to satisfy:

$$\left\| \Delta\overline{X}^{\,0p} \right\| \le \left(\beta_{pn}^{c}\right)^{-1}\left(\alpha_{pn}-\beta_{pn}^{c}\right)\left\| X-\overline{X}^{\,0p} \right\| \;,\; \forall X \in cl\,S \;,\; \forall n \in \boldsymbol{N} \tag{17}$$

Assume also that the nominal and incremental controller are

$$g_n\left(X_{n-1}\right)=\lambda_n\left(X_{n-1}\right)f_{n-\sigma_n}\left(X_{n-1-\sigma_n}\right) \;;\quad \tilde{g}_n\left(X_{n-1}\right)=\tilde{\lambda}_n\left(X_{n-1}\right)\tilde{f}_{n-\tilde{\sigma}_n}\left(X_{n-1-\tilde{\sigma}_n}\right)\;;\; \forall n \in \boldsymbol{N} \tag{18}$$

with their corresponding gain sequences $\left\{\lambda_n\right\}_{n\in\boldsymbol{N}}$ and $\left\{\tilde{\lambda}_n\right\}_{n\in\boldsymbol{N}}$ being chosen to satisfy the constraints:

$$sign\,\lambda_n\left(X_{n-1}\right)=-sign\left(\left(f_n\left(X_{n-1}\right)-\overline{x}^{\,0p}\right)f_{n-\sigma_n}\left(X_{n-1-\sigma_n}\right)\right)\;;\; \forall n \in \boldsymbol{N} \tag{19}$$

$$sign\,\tilde{\lambda}_n\left(X_{n-1}\right)=-sign\left(\left(\tilde{f}_n\left(X_{n-1}\right)-\overline{x}^{\,0p}\right)\tilde{f}_{n-\tilde{\sigma}_n}\left(X_{n-1-\tilde{\sigma}_n}\right)\right)\;;\; \forall n \in \boldsymbol{N} \tag{20}$$

$$\left\| \lambda_n\left(X_{n-1}\right), \tilde{\lambda}_n\left(X_{n-1}\right) \right\| \le \frac{\left| f_n\left(X_{n-1}\right)+\tilde{f}_n\left(X_{n-1}\right)-\overline{x}^{\,cp} \right|}{\left\| \left( f_{n-\sigma_n}\left(X_{n-1-\sigma_n}\right), \tilde{f}_{n-\tilde{\sigma}_n}\left(X_{n-1-\tilde{\sigma}_n}\right)\right) \right\|}\;;\; \forall n \in \boldsymbol{N} \tag{21}$$

for some existing non-negative integer sequences $\left\{\sigma_n\right\}_{n\in\boldsymbol{N}_0}$, $\left\{\tilde{\sigma}_n\right\}_{n\in\boldsymbol{N}_0}$ chosen such that $\left\| \left( f_{n-\sigma_n}\left(X_{n-1-\sigma_n}\right), \tilde{f}_{n-\tilde{\sigma}_n}\left(X_{n-1-\tilde{\sigma}_n}\right)\right) \right\| \ne 0$; $\forall n \in \boldsymbol{N}_0$ subject to $\liminf_{n\to\infty} \left\| \left( f_{n-\sigma_n}\left(X_{n-1-\sigma_n}\right), \tilde{f}_{n-\tilde{\sigma}_n}\left(X_{n-1-\tilde{\sigma}_n}\right)\right) \right\| \ge 0$ since equality to zero holds for all nonnegative sequences $\left\{\sigma_n\right\}_{n\in\boldsymbol{N}_0}$, $\left\{\tilde{\sigma}_n\right\}_{n\in\boldsymbol{N}_0}$ if and only if $\overline{X}^{\,cp}=0$.

Then, the corresponding equilibrium points of the perturbed uncontrolled associate vector system $\overline{X}^{\,0p}$ is unstable while that of the perturbed controlled system $\overline{X}^{\,cp}$ is asymptotically stable. The properties hold for the corresponding perturbed and perturbed controlled difference equations of equilibrium points $\overline{x}^{\,0p}$ and $\overline{x}^{\,cp}$, respectively .

**Proof**: Conditions for the following chain of inequalities to hold are given:

$$\begin{aligned} &\left| h_n\left(X_{n-1}\right)-\overline{x}^{\,cp} \right| \\ &\quad = \left| f_n\left(X_{n-1}\right)+\tilde{f}_n\left(X_{n-1}\right)+g_n\left(X_{n-1}\right)+\tilde{g}_n\left(X_{n-1}\right)-\overline{x}^{\,0p}-\Delta\overline{x}^{\,0p} \right| \\ &\quad \le \beta_{pn}^{c}\left\| X_{n-1}-\overline{X}^{\,c} \right\| \le \alpha_{pn}\left\| X_{n-1}-\overline{X}^{\,0p} \right\| \le \left| f_n\left(X_{n-1}\right)+\tilde{f}_n\left(X_{n-1}\right)-\overline{x}^{\,0p} \right| \end{aligned} \tag{22}$$

; $\forall n \in \boldsymbol{N}$ within $S$ .The following chained inequalities guarantee that (22) holds in $S$:

$$\begin{aligned} &\alpha_{pn}\left\| X-\overline{X}^{\,0p} \right\| - \left| \lambda_n\left(X_{n-1}\right)f_{n-\sigma_n}\left(X_{n-1-\sigma_n}\right)+\tilde{\lambda}_n\left(X_{n-1}\right)\tilde{f}_{n-\tilde{\sigma}_n}\left(X_{n-1-\tilde{\sigma}_n}\right)-\Delta\overline{x}^{\,0p} \right| \\ &\quad \le \left| f_n\left(X_{n-1}\right)+\tilde{f}_n\left(X_{n-1}\right)-\overline{x}^{\,0p}+\lambda_n\left(X_{n-1}\right)f_{n-\sigma_n}\left(X_{n-1-\sigma_n}\right)+\tilde{\lambda}_n\left(X_{n-1}\right)\tilde{f}_{n-\tilde{\sigma}_n}\left(X_{n-1-\tilde{\sigma}_n}\right)-\Delta\overline{x}^{\,0p} \right| \\ &\quad = \left| f_n\left(X_{n-1}\right)+\tilde{f}_n\left(X_{n-1}\right)-\overline{x}^{\,0p}-\Delta\overline{x}^{\,0p} \right| - \left| \lambda_n\left(X_{n-1}\right)f_{n-\sigma_n}\left(X_{n-1-\sigma_n}\right)+\tilde{\lambda}_n\left(X_{n-1}\right)\tilde{f}_{n-\tilde{\sigma}_n}\left(X_{n-1-\tilde{\sigma}_n}\right) \right| \\ &\qquad \le \beta_{pn}^{c}\left\| X-\overline{X}^{\,0p}-\Delta\overline{X}^{\,0p} \right\| \le \alpha_{pn}\left\| X-\overline{X}^{\,0p} \right\|;\; \forall n \in \boldsymbol{N} \end{aligned} \tag{23}$$

$$0 \le \left| f_n(X_{n-1}) + \tilde{f}_n(X_{n-1}) - \bar{x}^{0p} - \Delta\bar{x}^{0p} \right| - \beta_{pn}^c \left\| X - \bar{X}^{0p} - \Delta\bar{X}^{0p} \right\|$$
$$\le \left| \lambda_n(X_{n-1}) f_{n-\sigma_n}(X_{n-1-\sigma_n}) + \tilde{\lambda}_n(X_{n-1}) \tilde{f}_{n-\tilde{\sigma}_n}(X_{n-1-\tilde{\sigma}_n}) \right| \; ; \; \forall n \in \boldsymbol{N} \qquad (24)$$

Then, the nominal and incremental controller gains are chosen to satisfy (19)-(21) for existing non-negative real sequences $\{\sigma_n\}_{n\in \boldsymbol{N}_0}$ , $\{\tilde{\sigma}_n\}_{n\in \boldsymbol{N}_0}$ such that $\left\| \left( f_{n-\sigma_n}(X_{n-1-\sigma_n}),\, \tilde{f}_{n-\tilde{\sigma}_n}(X_{n-1-\tilde{\sigma}_n}) \right) \right\| \ne 0$ ; $\forall n \in \boldsymbol{N}_0$ subject to $\liminf_{n\to\infty} \left\| \left( f_{n-\sigma_n}(X_{n-1-\sigma_n}),\, \tilde{f}_{n-\tilde{\sigma}_n}(X_{n-1-\tilde{\sigma}_n}) \right) \right\| \ge 0$ since equality to zero holds if and only if $\bar{X}^{cp} = 0$ . The reminder of the proof now follows since from (22) one gets simultaneously within a nonempty invariant set $S \subseteq A_{0p} \cap A_{cp}$ :

$$\left| h_n(X_{n-1}) - \bar{x}^{cp} \right| \le \beta_{pn}^c \left\| X_{n-1} - \bar{X}^{cp} \right\| ; \; \forall X_i \in S \; , \; \forall i \in \boldsymbol{N}_0 \qquad (25)$$

$$\alpha_{pn} \left\| X_{n-1} - \bar{X}^{0p} \right\| \le \left| f_n(X_{n-1}) + \tilde{f}_n(X_{n-1}) - \bar{x}^{0p} \right| ; \; \forall X_i \in S \; , \; \forall i \in \boldsymbol{N}_0 \qquad (26)$$

Then, one gets from (26):

$$\left| h_n(X) - \bar{x}^{cp} \right| = \left| h_n\left(V_{h_{n-1}}(X)\right) - \bar{X}^{cp} \right| = \left| h_n\left(G_{n-1}^h(X)\right) - \bar{x}^{cp} \right|$$
$$\le \beta_{pn}^c \left\| X - \bar{X}^{cp} \right\| \le \prod_{i=1}^n \left[ \max_{i \ge j(\in \boldsymbol{N}) \ge i-m-1} \left( \beta_{pj}^c \right) \right] \left\| X - \bar{X}^{cp} \right\| \le \bar{\beta}^{cn} \left\| X - \bar{X}^{cp} \right\| ; \forall X \in cl\, S \; , \forall n \in \boldsymbol{N} \qquad (27)$$

with the real sequence $\left\{ \bar{\beta}^{cn} \right\}_{n \in \boldsymbol{N}_0}$ of elements satisfying $\bar{\beta}^{cn} \in [0,1)$, $\forall n \in \boldsymbol{N}$ so that one deduces by taking $\ell_\infty$ -norms for the m-tuples $G_{(.)}^h(X)$ (see Remark 2.2) that:

$$\left\| G_{n-1}^h(X) - \bar{X}^{cp} \right\|_\infty = \max_{n \ge j(\in \boldsymbol{N}) \ge n-m-1} \left| h_j\left( G_{j-1}^h(X) \right) - \bar{x}^{cp} \right| \le \bar{\beta}^{cn} \left\| X - \bar{X}^{cp} \right\|_\infty ; \forall n \in \boldsymbol{N} \qquad (28)$$

Since $\bar{\beta}^{cn} \in (0,1)$; $\forall n \in \boldsymbol{N}$ , then $\lim_{n\to\infty} \left( G_{n-1}^h(X) - \bar{X}^{cp} \right) = 0$ ; $\forall X \in S$ . Since $\left\{ G_n^h(X) \right\}_{n \in \boldsymbol{N}}$ is a sequence of contraction self-mappings from $\boldsymbol{R}^m \,|\, S$ to S and $\boldsymbol{R}^m$ is a complete metric space endowed with the given norm-induced metric, then the equilibrium point $\bar{X}^{cp}$ on the vector function equation associated to the controlled difference equation is locally asymptotically stable with respect to S and it also is the unique fixed point in $cl\, S$ ,[8]. Then, the equilibrium point $\bar{X}^{cp}$ on the controlled difference equation is locally asymptotically stable with respect to S equilibrium point $\bar{x}^{cp}$ of the corresponding difference equation is also locally stable. On the other hand, it follows from (26) that

$$\alpha_{pn-1} \left\| X - \bar{X}^{0p} \right\|_\infty = \alpha_{pn} \left\| G_{n-1}^{f+\tilde{f}}(X) - \bar{X}^{0p} \right\|_\infty$$
$$\le \left| f_n(X) + \tilde{f}_n(X) - \bar{x}^{0p} \right| = \left| f_n\left(V_{f_n+\tilde{f}_n}(X)\right) + \tilde{f}_n\left(V_{f_n+\tilde{f}_n}(X)\right) - \bar{x}^{0p} \right| \le \left\| G_n^{f+\tilde{f}}(X) - \bar{X}^{0p} \right\|_\infty ; \forall X \in cl\, S \; , \forall n \in \boldsymbol{N} \qquad (29)$$

Since $\{\overline{\alpha}_{pn}\}_{n\in\mathbf{N}}$ is unbounded with $\overline{\alpha}_{pn}:=\prod_{i=1}^{n}\alpha_{pi}$, then $\left\{G_n^{f+\tilde{f}}(X)\right\}_{n\in\mathbf{N}}$ is a sequence of expanding self-mappings from $\mathbf{R}^m|S$ to S so that the equilibrium point $\overline{X}^{0p}$ of the vector function equation associated to the perturbed uncontrolled difference equation is locally unstable with respect to S. Then, the corresponding equilibrium point of the difference equation $\overline{x}^{0p}$ is unstable. □

The next result is concerned with the local asymptotic stabilization around equilibrium points within a certain invariant set of the unstable uncontrolled perturbed difference equation through a single controller provided that the perturbation dynamics satisfies some smallness- type constraints.

**Theorem 3.3**. Assume the following:

1) $m=m_{0p}=m_c$ with $\overline{X}^{0p}\in cl\,S\cap A_{0p}$ and $\overline{X}^{c}\in cl\,S\cap A_c$ being unique equilibrium points in $cl\,S\cap A_{0p}$, respectively in $cl\,S\cap A_c$, where $S\subseteq A_{0p}\cap A_c$ is invariant under all $V_{f_n+\tilde{f}_n}$ and $V_{h_n}\equiv V_{f_n+\tilde{f}_n+g_n}$ (since $\tilde{g}_n\equiv 0$) for any controller in the class $\boldsymbol{C}$ consisting of a nominal controller; i.e. $V_{f_n+\tilde{f}_n}(S)\subseteq S$ and $V_{h_n}(S)\subseteq S$ with the incremental controller being identically zero; $\forall n\in\mathbf{N}$, provided that the sets are nonempty. Assume also that the sets $A_{0p}$ and $A_c$ defined in (14)-(15) are redefined as

$$A_{0p}:=\left\{X\in R^m:\left|\,f_n(X)+\tilde{f}_n(X)-\overline{x}^{0p}\right|\geq\alpha_{pn}\left(1-\tilde{\alpha}_{pn}\right)\left\|X-\overline{X}^{0p}\right\|;\forall n\in\mathbf{N}\right\}\tag{30}$$

$$A_c:=\left\{X\in R^m:\left|\,f_n(X)+\tilde{f}_n(X)+g_n(X)-\overline{x}^{c}\right|\leq\beta_n^c\left(1+\tilde{\beta}_n^c\right)\left\|X-\overline{X}^{c}\right\|;\forall n\in\mathbf{N}\right\}\tag{31}$$

for nonnegative real sequences $\{\alpha_{pn}(1-\tilde{\alpha}_{pn})\}_{n\in\mathbf{N}}$ and $\left\{\beta_n^c\left(1+\tilde{\beta}_n^c\right)\right\}_{n\in\mathbf{N}}$ being defined for some nonnegative real sequences $\{\tilde{\alpha}_{pn}\}_{n\in\mathbf{N}}$ and $\left\{\tilde{\beta}_n^c\right\}_{n\in\mathbf{N}}$, subject to $\tilde{\alpha}_{pn}\leq 1$ and $\tilde{\beta}_n^c<\beta_n^{-1}-1$, where $\{\overline{\alpha}_{pn}\}_{n\in\mathbf{N}}$ is unbounded of elements redefined as $\overline{\alpha}_{pn}:=\prod_{i=1}^{n}\left[\alpha_{pi}\left(1-\tilde{\alpha}_{pi}\right)\right]$, and $\{\overline{\beta}_n^c\}_{n\in\mathbf{N}}$ of elements being redefined as $\overline{\beta}_n^c:=\prod_{i=1}^{n}\left[\beta_i^c\left(1+\tilde{\beta}_i^c\right)\right]$ being such that $\limsup_{n\to\infty}\overline{\beta}_n^c<1$.

2) The perturbed sequence $\left\{\tilde{f}_n(X_{n-1})\right\}_{n\in\mathbf{N}}$ satisfies the constraints

$$\left|\,\tilde{f}_n(X_{n-1})-a_n\right|\leq\beta_n^c\tilde{\beta}_n^c\left\|X_{n-1}-\overline{X}^{c}\right\|\;;\;\left|\,\tilde{f}_n(X_{n-1})-b_n\right|\leq\alpha_{pn}\tilde{\alpha}_{pn}\left\|X_{n-1}-\overline{X}^{0p}\right\|\tag{32}$$

; $\forall n\in\mathbf{N}$ within $S$ for some real sequences $\{a_n\}_{n\in\mathbf{N}}$ and $\{b_n\}_{n\in\mathbf{N}}$ and for some nonnegative real sequences $\{\tilde{\alpha}_{pn}\}_{n\in\mathbf{N}_0}$ and $\{\tilde{\beta}_n\}_{n\in\mathbf{N}_0}$, subject to $\tilde{\alpha}_{pn}\leq 1$ and $\tilde{\beta}_{pn}<\beta_n^{-1}-1$.

3) The stabilizing incremental controller is identically zero while the nominal controller is $g_n(X_{n-1})=\lambda_n(X_{n-1})f_{n-\sigma_n}(X_{n-1-\sigma_n})$ subject to:

$$sign\,\lambda_n(X_{n-1}) = -sign\left(\left(f_n(X_{n-1}) - \bar{x}^{0p}\right) f_{n-\sigma_n}\left(X_{n-1-\sigma_n}\right)\right); \forall n \in \boldsymbol{N} \tag{33}$$

$$\left|\lambda_n(X_{n-1})\right| \le \frac{\left|f_n(X_{n-1}) - \bar{x}^{c} + a_n\right|}{\left|f_{n-\sigma_n}\left(X_{n-1-\sigma_n}\right)\right|}; \ \forall n \in \boldsymbol{N} \tag{34}$$

within S if $f_{n-\sigma_n}\left(X_{n-1-\sigma_n}\right) \neq 0$ and $\lambda_n(X_{n-1}) = 0$ if $f_{n-\sigma_n}\left(X_{n-1-\sigma_n}\right) = 0$; $\forall n \in \boldsymbol{N}$, for some existing non-negative integer sequence $\{\sigma_n\}_{n \in \boldsymbol{N}_0}$.

Then, the corresponding equilibrium points of the perturbed uncontrolled associate vector system $\bar{X}^{0p}$ is unstable while that of the controlled system $\bar{X}^{c}$ under the nominal controller is asymptotically stable. The stability properties also hold for the corresponding perturbed and perturbed controlled difference equations of equilibrium points $\bar{x}^{0p}$ and $\bar{x}^{c}$, respectively .

**Proof**: Now the errors of the uncontrolled perturbed and the controlled nominal equilibrium points under consideration are $\Delta\bar{x}^{0p} := \bar{x}^{c} - \bar{x}^{0p}$ and $\Delta\bar{X}^{0p} := \bar{X}^{c} - \bar{X}^{0p}$. Note that $\tilde{g}_n(X_{n-1}) \equiv 0$, so that (20) is omitted, while (19) and (21) are replaced by (33) and (34).Also, the inequalities (25)-(26) are replaced by

$$\left| f_n(X_{n-1}) + \tilde{f}_n(X_{n-1}) + \lambda_n(X_{n-1}) f_{n-\sigma_n}\left(X_{n-1-\sigma_n}\right) - \bar{x}^{c} \right| \le \left| f_n(X_{n-1}) + \lambda_n(X_{n-1}) f_{n-\sigma_n}\left(X_{n-1-\sigma_n}\right) - \bar{x}^{c} + a_n \right|$$
$$+ \left| \tilde{f}_n(X_{n-1}) - a_n \right| \le \beta_n^c\left(1 + \tilde{\beta}_n^c\right)\left\| X_{n-1} - \bar{X}^{cp} \right\|; \ \forall X_n \in S, \forall n \in \boldsymbol{N}_0 \tag{35}$$

$$\alpha_{pn}\left(1 - \tilde{\alpha}_{pn}\right)\left\| X_{n-1} - \bar{X}^{0p} \right\| \le \left| f_n(X_{n-1}) - \bar{x}^{0p} + b_n \right| - \left| \tilde{f}_n(X_{n-1}) - b_n \right|$$
$$\le \left| f_n(X_{n-1}) + \tilde{f}_n(X_{n-1}) - \bar{x}^{0p} \right|; \ \forall X_n \in S, \ \forall n \in \boldsymbol{N}_0 \tag{36}$$

$\forall n \in \boldsymbol{N}_0$, provided that

$$\left| f_n(X_{n-1}) + \lambda_n(X_{n-1}) f_{n-\sigma_n}\left(X_{n-1-\sigma_n}\right) - \bar{x}^{c} + a_n \right| \le \beta_n^c \left\| X_{n-1} - \bar{X}^{cp} \right\|; \ \forall X_n \in S, \forall n \in \boldsymbol{N}_0 \tag{37}$$

$$\alpha_{pn}\left\| X_{n-1} - \bar{X}^{0p} \right\| \le \left| f_n(X_{n-1}) - \bar{x}^{0p} + b_n \right|; \ \forall X_n \in S, \forall n \in \boldsymbol{N}_0 \tag{38}$$

fo some existing nonnegative sequences $\left\{\alpha_{pn}\left(1 - \tilde{\alpha}_{pn}\right)\right\}_{n \in \boldsymbol{N}}$ and $\left\{\beta_n^c\left(1 + \tilde{\beta}_n^c\right)\right\}_{n \in \boldsymbol{N}}$. Thus, $\bar{X}^{0p}$ is unstable while $\bar{X}^{c}$ is locally asymptotically stable with respect to *S*. Those properties also hold by construction of $\bar{X}^{0p}$ and $\bar{X}^{c}$ for the corresponding perturbed and perturbed controlled difference equations of equilibrium points $\bar{x}^{0p}$ and $\bar{x}^{c}$, respectively . □

*Remark 3.1.* Note that in the proof of the results of this section, estimates can replace to the true equilibrium points if they are , potentially distinct, but sufficiently close to each other by using the results of Section 2 provided that the needed assumptions of the various function smoothness hold. Furthermore, the equilibrium points under analysis in the various given results could be replaced with the estimates of errors related to the nominal equilibrium if such errors are sufficiently small in terms of smallness of error norms. For instance, take the estimation error of the equilibrium points of the feedback associate vector equation via the nominal controller compared to its uncontrolled perturbed counterpart

$$\Delta := \hat{\tilde{X}}^{c} - X^{0} + \left(I_m - M^{c}\left(X^{0}\right)\right)^{-1} V_{\tilde{f}+g}\left(X^{0}\right)$$

; $\forall X \in S$ calculated from Theorem 2.6. Thus, the last term of (35) with $\tilde{\beta}_n^c = 0$ ; $\forall n \in \boldsymbol{N}_0$ possesses the lower- bounds given below

$$\beta_n^c \left\| X - \overline{X}^c \right\| = \beta_n^c \left\| X - \hat{\tilde{X}}^c - \Delta \right\| \geq \beta_n^c \left\| X - \hat{\tilde{X}}^c \right\| - \beta_n^c \left\| \Delta \right\| \geq \beta_n^c \left(1 - \varepsilon_n^c\right) \left\| X - \hat{\tilde{X}}^c \right\| \qquad (39)$$

provided that the equilibrium point estimation error is of sufficiently small size fulfilling $\left\| \Delta \right\| \leq \varepsilon_n^c \left\| X - \hat{\tilde{X}}^c \right\|$ for some sequence $\left\{ \varepsilon_n^c \right\}_{n \in \boldsymbol{N}_0}$ satisfying $0 \leq \varepsilon_n^c \leq 1$ ; $\forall n \in \boldsymbol{N}_0$. Thus, the last term of (35) could be replaced by (39) and the theorem and its proof could be reformulated based on estimates when having a sufficiently small estimation error between the equilibrium point of the uncontrolled difference equation and that of the current controlled difference equation; i.e. that subject to parametrical perturbation with or without unmodeled dynamics. □

**ACKNOWLEDGMENTS**

The author thanks to the Spanish Ministry of Education by its support of this work through Grant DPI2009-07197. He is also grateful to the Basque Government by its support through Grant GIC07143-IT-269-07.

**REFERENCES**

[1] R. Memarbashi, "On the behaviour of nonautonomous difference equations", *Journal of Difference Equations and Applications*, Vol. 16, Issue 9, pp. 1031-1035, 2010.

[2] R. Memarbashi, "Sufficient conditions for the exponential stability of nonautonomous difference equations", *Applied Mathematics Letters*, Vol. 21, pp. 232-235, 2008.

[3] R. Memarbashi, "On the stability of nonautonomous difference equations", *Journal of Difference Equations and Applications*, Vol. 14, No. 3, pp. 301-307, 2008.

[4] E. Liz and J.B. Ferreiro, "A note on the global stability of generalised difference equations", *Applied Mathematics Letters* , Vol. 15, pp. 755-659, 2002.

[5] H. Sedaghat, "Geometric stability conditions for higher order difference equations", *Journal of Mathematical Analysis with Applications*, Vol. 224, pp. 255-272, 1998.

[6] R.P. Agarwal, *Difference Equations and Inequalities*, 2nd ed., Dekker, New York, 2000.

[7] M. Furi, M. Martelli and M. O´Neill, "Global stability of equilibria", *Journal of Difference Equations and Applications*, Vol. 15, issue 4, pp. 387-397, 2009.

[8] M. De la Sen, "Stable iteration procedures in metric spaces which generalize a Picard-type iteration", *Fixed Point Theory and Applications*, Article ID 953091, 2010.

[9] S. Alonso-Quesada and M. De la Sen, "Robust adaptive control of discrete nominally stabilizable plants", *Applied Mathematics and Computation*, Vol. 150, Issue 2, pp. 555-583. 2004.

[10] Q.M. Liu, R. Xu and Z.P. Wang, "Stability and bifurcation of a class of discrete-time Cohen-Grossberg neural networks with delays", *Discrete Dynamics in Nature and Society*, Article ID 403873, 2011.

[11] S.X. Ye and W.Q. Wang, "Stability analysis and stabilisation for a class of 2-D nonlinear discrete systems", *International Journal of Systems Science*, Vol. 42, Issue 5, pp. 839-851, 2011.

[12] A. Dey and H. Kar, "Robust stability of 2-D discrete systems employing generalized overflow nonlinearities: An LMI approach", *Digital Signal Processing*, Vol. 21, Issue 2, pp. 187-195, 2011.

[13] E.B.M. Bashier and K.C. Patidar, "A novel fitted operator finite difference method for a singularly perturbed delay parabolic partial difference equation", *Applied Mathematics and Computation* , Vol. 217, Issue 9, pp. 4728-4739, 2011.

[14] U. Krause, "Stability of non-autonomous population models with bounded and periodic enforcement", *Journal of Difference Equations and Applications*, Vol.15, Issue 7, pp. 649-658, 2009.

[15] C.Y. Wang, S. Wang and W. Wang, " Global asymptotic stability of equilibrium point for a family of rational difference equations", *Applied Mathematics Letters*, Vol. 24, Issue 5, pp. 714-718, 2011.

[16] J.M. Ortega, *Numerical Analysis. A Second Course*, Classics in Applied Mathematics, Vol. 3, SIAM, Philadelphia, PA, 1990.

[17] S. Elaydi, *An Introduction to Difference Equations*, 2nd ed., Springer-Verlag, New York, 1999.

[18] M.R.S. Kulenovic and G. Ladas, *Dynamics of Second Order Rational Difference Equations with Open Problems and Conjectures*, Chapman and Hall/CRC Press, London/ Boca Raton, FL., 2002.